\documentclass[reqno,a4paper]{amsart}

\usepackage{amsmath,amsfonts,amssymb}
\usepackage{verbatim}
\usepackage{enumerate}
\usepackage{url}
\usepackage{tikz}
\usepackage{setspace}
\usepackage[left=2.5cm,right=2.5cm, top=3cm, bottom=2.5cm]{geometry}
\usepackage[utf8]{inputenc} 
\usepackage[T1]{fontenc}
\usepackage{hyperref} 	
\usetikzlibrary{arrows}

\def\d{{\sf d}}

\def\s{{\sf s}}
\def\e{{\sf \eta}}
\def\E{{\sf E}}
\def\F{{\mathcal F}}
\def\bd{{\boldsymbol{\cdot}}}

\def\Z{\mathbb Z}

\makeatletter
\newcommand{\dotprod}{\mathbin{\mathpalette\dotprod@\relax}}
\newcommand{\dotprod@}[2]{%
  \ooalign{$\m@th#1{\LARGE\bd}$\cr\hidewidth$\m@th#1\prod$\hidewidth\cr}%
}
\def\ldotprod{\mathop{\dotprod}\limits}
\makeatother

\theoremstyle{plain}
\newtheorem{theorem}{Theorem}[section]
\newtheorem{lemma}[theorem]{Lemma}

\newtheorem{conjecture}[theorem]{Conjecture}
\def\proof{\noindent {\it Proof: }}
\def\qed{\hfill\hbox{$\square$}}

\theoremstyle{definition}

\numberwithin{equation}{section}

\title{On the direct and inverse zero-sum problems over $C_n \rtimes_s C_2$}

\author[D.V. Avelar]{D.V. Avelar}
\address{Departamento de Análise\\
Universidade Federal Fluminense (UFF)\\
Niterói, RJ\\
24210-201\\
Brazil\\
}
\email{daniloavelar@id.uff.br}

\author[F.E. Brochero Mart\'{\i}nez]{F.E. Brochero Mart\'{\i}nez$^*$}
\address{
Departamento de Matem\'{a}tica\\
Universidade Federal de Minas Gerais (UFMG)\\
Belo Horizonte, MG\\
31270-901\\
Brazil\\
}
\email{fbrocher@mat.ufmg.br}

\author[S. Ribas]{S. Ribas$^*$}
\address{
Departamento de Matem\'{a}tica\\
Universidade Federal de Ouro Preto (UFOP)\\
Ouro Preto, MG\\
35400-000\\
Brazil\\
}
\email{savio.ribas@ufop.edu.br}

\thanks{$^*$Partially supported by FAPEMIG APQ-02546-21, Brazil.}

\date{\today}

\onehalfspace

\usepackage{xcolor}

\begin{document}

\maketitle

\begin{abstract}
Let $C_n$ be the cyclic group of order $n$. In this paper, we provide the exact values of some zero-sum constants over $C_n \rtimes_s C_2$ where $s \not\equiv \pm1 \pmod n$, namely $\e$-constant, Gao constant, and Erd\H os-Ginzburg-Ziv constant (the latter for all but a ``small'' family of cases). As a consequence, we prove the Gao's and Zhuang-Gao's Conjectures for groups of this form. We also solve the associated inverse problems by characterizing the structure of product-one free sequences over $C_n \rtimes_s C_2$ of maximum length.
\end{abstract}

\section{Introduction}

For a finite multiplicative group $G$, the {\em zero-sum problems} consist of stablishing conditions that guarantee that a given sequence over $G$ has a non-empty product-one subsequence with some prescribed property. This kind of problem dates back to the celebrated works of Erd\H os, Ginzburg \& Ziv \cite{EGZ}, van Emde Boas \& Kruyswijk \cite{vEBK} and Olson \cite{Ols1,Ols2}. It has applications and connections in several branches of mathematics; we refer to the surveys from Caro \cite{Car} and Gao \& Geroldinger \cite{GaGe} for an overview over abelian groups. Decades later, these kind of problems were further generalized to non-abelian groups with the works \cite{Yu,YP,ZhGa}. This explains the multiplicative notation; in particular, we use ``product-one'' rather than ``zero-sum''.

\subsection{Definitions and notations}

Let $\F(G)$ be a free abelian monoid, written multiplicatively, with basis $G$. A sequence $S$ over $G$ is a finite and unordered element of $\F(G)$, which is equipped with the sequence concatenation product denoted by ${\bd}$. Therefore, a sequence $S \in \F(G)$ has the form 
$$S = g_1 {\bd} \dots {\bd} g_k = \ldotprod_{1\le i \le k}
\, g_i = \ldotprod_{1\le i \le k}
g_{\tau(i)},$$ 
for any permutation $\tau: \{1,2,\dots,k\} \to \{1,2,\dots,k\}$, where $g_1, \dots, g_k \in G$ are the {\em terms} of $S$ and $k = |S| \ge 0$ is the {\em length} of $S$. Given $g \in G$ and $t \ge 0$, we shorten $g^{[t]} = \ldotprod_{1\le i \le t} g$. 
The {\em multiplicity} of the term $g \mid S$ is denoted by $v_g(S) = \#\{i \in \{1,2,\dots,k\} \mid g_i = g\}$, thus we may write $S = \ldotprod_{g \in G} g^{[v_g(S)]}$. 
A sequence $T$ is a {\em subsequence} of $S$ if $T \mid S$ as elements of $\F(G)$; equivalently, $v_g(T) \le v_g(S)$ for all $g \in G$. In this case, we write $S {\bd} T^{[-1]} = \ldotprod_{g \in G}  g^{[v_g(S) - v_g(T)]}$ and $T^{[k]} = \ldotprod_{1\le i \le k}T$. Moreover, let $S \cap K = \ldotprod_{g \mid S \atop g \in K} g^{[v_g(S)]}$ be the subsequence of $S$ formed by the terms that lie in a subset $K \subset G$.

Let $H$ be a normal subgroup of $G$ and let $\overline{G} = G/H$. Consider the natural homomorphism $\phi: G \to \overline{G}$. If $g \in G$, we denote $\overline{g} = \phi(g)$. In this way, we write $\overline{S} = \ldotprod_{1 \le i \le k} \overline{g_i} \in \F(\overline{G})$.

Denote the {\em set of products of $S$}, {\em the set of subproducts of $S$}, and {\em the set of $n$-subproducts of $S$}, respectively, by
$$\pi(S) = \left\{ \prod_{i=1}^k g_{\tau(i)} \in G \mid \tau \text{ is a permutation of $[1,k]$}\right\}, \quad \Pi(S) = \bigcup_{T \mid S \atop |T| \ge 1} \pi(T), \quad \text{and} \quad \Pi_n(S) = \bigcup_{T \mid S \atop |T| = n} \pi(T).$$


The sequence $S$ is:
\begin{itemize}
\item {\em product-one sequence} if $1 \in \pi(S)$,
\item {\em short product-one sequence} if $S$ is a product-one sequence with $1 \le |S| \le \exp(G)$,
\item {\em product-one free} if $1 \not\in \Pi(S)$,
\item {\em $n$-product-one free} if $1 \not\in \Pi_n(S)$.
\end{itemize}

Furthermore, the following zero-sum invariants are defined:
\begin{enumerate}[(a)]
\item {\em Small Davenport constant}, $\d(G)$, is the maximum length of a product-one free sequence over $G$,
\item {\em Gao constant}, $\E(G)$, is the smallest $\ell > 0$ such that every $S \in \F(G)$ with $|S| \ge \ell$ has a $|G|$-product-one subsequence,
\item {\em $\e$-constant}, $\e(G)$, is the smallest $\ell > 0$ such that every $S \in \F(G)$ with $|S| \ge \ell$ has a short product-one subsequence,
\item {\em Erd\H os-Ginzburg-Ziv constant}, $\s(G)$, is the smallest $\ell > 0$ such that every $S \in \F(G)$ with $|S| \ge \ell$ has an $\exp(G)$-product-one subsequence.
\end{enumerate}



\subsection{Background}

The previous invariants are well-defined and finite for finite groups $G$. Indeed, Pigeonhole Principle ensures that $\d(G) +1 \le 
|G|$. 
This inequality can be further improved for non-abelian groups: $\d(G) \le \lfloor |G|/2 \rfloor$ \cite{OW}. 
By definition, it follows that $\d(G) \le \e(G) \le \s(G) \le \E(G)$. Moreover $\E(G) \le 2|G|-1$ \cite[Theorem~10.1]{Gr1}. Other non-trivial bounds (see \cite[Lemma 4]{ZhGa}) are 
\begin{equation}\label{lowerboundse}
\s(G) \ge \e(G)+\exp(G)-1 \quad \text{ and } \quad \E(G) \ge \d(G) + |G|.
\end{equation}

\begin{conjecture}[Gao \cite{Ga2}]\label{conjgao}
$\s(G) = \e(G) + \exp(G) - 1$ for every finite group $G$.
\end{conjecture}

\begin{conjecture}[Zhuang-Gao \cite{ZhGa}]\label{conjzhuanggao}
$\E(G) = \d(G) + |G|$ for every finite group $G$.
\end{conjecture}

In \cite{Ga}, Conjecture \ref{conjzhuanggao} was proven for abelian groups. After finding the exact values of those invariants, which is called {\emph{direct problem}}, it is natural to consider the related problem of finding the structure of the product-one free sequences of maximum length with some prescribed property. This is called {\emph{inverse problem}}. In \cite[Theorem 3.2]{Sc} and \cite[Theorem 11.1]{Gr1}, the inverse problems related to small Davenport constant for abelian groups of ranks one and two are solved, respectively. In \cite{Ga1} (see Lemma \ref{lemmaegzinversecyclic}), it is solved the inverse problem related to Erd\H os-Ginzburg-Ziv constant for cyclic groups, and this result will be widely used thoughout this paper. For non-abelian groups, we refer to \cite{MR2,OhZh,MR1,MR3,QL} for some recent results in this direction.

Let $C_n$ be the cyclic group of order $n$, and $D_{2n}$ be the dihedral group of order $2n$. The values of $\d(G)$, $\e(G)$, $\s(G)$ and $\E(G)$ are well-known for abelian groups of rank at most two and for some classes of non-abelian groups. In particular, we have that 
\begin{itemize}
\item $\d(C_n) = n-1$, 
$\e(C_n) = n$ \cite{Ols1}, $\s(C_n) = \E(C_n) = 2n-1$ \cite{EGZ}.
\item 
$\d(D_{2n}) = n$ \cite{ZhGa}, $\e(D_{2n}) = n+1$ \cite{Zh}, $\s(D_{2n}) = \begin{cases} 2n \text{ if $n$ is even \cite{OhZh},} \\ 3n \text{ if $n$ is odd \cite{Bas},} \end{cases}$ and $\E(D_{2n}) = 3n$ \cite{Bas}.
\end{itemize}

%




\subsection{On the group $C_n \rtimes_s C_2$ and the main results}

Let 
\begin{equation}\label{Gns}
G_{n,s} = C_n \rtimes_s C_2 = \langle x,y \mid x^2 = y^n = 1, \; yx = xy^s \rangle.
\end{equation}
It is known that $G_{n,s}$ is indeed a group if and only if $s^2 \equiv 1 \pmod n$. Some remarkable groups of this form are, for instance, the abelian group $G_{n,1} \simeq C_n \times C_2$, the dihedral group $G_{n,-1}$, the modular maximal-cyclic group $G_{2^t, 2^{t-1}+1}$ (for $t \ge 3$) of order $2^{t+1}$, and the quasidihedral group $G_{2^t, 2^{t-1}-1}$ (for $t \ge 3$) of order $2^{t+1}$.
However there are several other groups of this form. We refer to \cite{Ols2,Sc,ZhGa,MR2,OhZh} as the main references for direct and inverse zero-sum problems over $G_{n,1}$ and over $G_{n,-1}$. For a general $G_{n,s}$, we refer to \cite{GeGr}, where the large Davenport constant was found, and to \cite{MR3}, where two of the authors solved the inverse problem related to the small Davenport constant.

In this paper, we deal with groups $G_{n,s}$ for which $s \not\equiv \pm1 \pmod n$. In view of Lemma \ref{lemmafactorn}, it is possible to factor $n = n_1n_2$ or $n = 2n_1n_2$, where $n_1, n_2$ are coprime integers with $n_1 \neq 2$, $n_2 \neq 2$, $s \equiv -1 \pmod {n_1}$ and $s \equiv 1 \pmod {n_2}$. In particular, we solve the direct and inverse problems related to the $\e$-constant, Erd\H os-Ginzburg-Ziv constant and Gao constant over $G_{n,s}$, except for the case where $n_1 = 3$ and $n_2$ is odd in the last two constants. The main results are stated below.

\begin{theorem}\label{thmeta}
Let $n \ge 8$ and $s$ be integers such that $s^2 \equiv 1 \pmod n$, but $s \not\equiv \pm1 \pmod n$. We have that $\eta(G_{n,s}) = n+1$. Moreover, let $S \in \F(G_{n,s})$ of length $|S| = n$. Then $S$ has no short product-one subsequences if and only if there exist $\alpha, \beta \in G_{n,s}$ such that $G_{n,s} \simeq \langle \alpha, \beta \mid \alpha^2 = \beta^n = 1, \beta \alpha = \alpha \beta^s \rangle$ and $S = \beta^{[n-1]} \bd (\alpha \beta^t)$ for some $0 \le t \le n-1$.
\end{theorem}

We highlight that these extremal sequences are precisely the same found in \cite{MR2} for the small Davenport constant (see also Lemma \ref{lemmadavinversegns}).

For $n$ odd, since $\exp(G_{n,s}) = 2n$, it follows that $\s(G_{n,s}) = \E(G_{n,s})$. On the other hand, for $n$ even, we have $\exp(G_{n,s}) = n$ and, by extracting two disjoint $n$-product-one subsequences, it follows that $\E(G_{n,s}) \le \s(G_{n,s}) + n$. Except for a small case, the precise values are given by the following result.

\begin{theorem}\label{thmesdirect}
Let $n \ge 8$ and $s$ be integers such that $s^2 \equiv 1 \pmod n$, but $s \not\equiv \pm1 \pmod n$. Suppose additionally that if $n$ is odd, then $n_1 > 3$, where $n_1$ is given by Lemma \ref{lemmafactorn}. We have that $\E(G_{n,s}) = 3n$, and $\s(G_{n,s}) = \begin{cases} 2n & \text{ if $n$ is even} \\ 3n & \text{ if $n$ is odd}\end{cases}$.
\end{theorem}

As immediate consequences, we obtain that Conjectures \ref{conjgao} and \ref{conjzhuanggao} hold true for $G_{n,s}$, except possibly for a small case. Furthermore, we also solve the associated inverse problems.

\begin{theorem}\label{thmesinverse}
Let $n \ge 8$ and $s$ be integers such that $s^2 \equiv 1 \pmod n$, but $s \not\equiv \pm1 \pmod n$. Consider $S \in \F(G_{n,s})$.
\begin{enumerate}[(a)]
\item Let $n$ be even, and $|S| = \s(G_{n,s}) - 1 = 2n-1$. Then $S$ is $n$-product-one free if and only if there exist $\alpha, \beta \in G_{n,s}$, $t_1, t_2, t_3 \in \Z$ such that $G_{n,s} \simeq \langle \alpha, \beta \mid \alpha^2 = \beta^n = 1, \beta \alpha = \alpha \beta^s \rangle$, $\gcd(t_1-t_2,n) = 1$ and $S = (\beta^{t_1})^{[n-1]} \bd (\beta^{t_2})^{[n-1]} \bd (\alpha \beta^{t_3})$.

\item Let $n$ be even, and $|S| = \E(G_{n,s}) - 1 = 3n-1$. Then $S$ is $2n$-product-one free if and only if there exist $\alpha, \beta \in G_{n,s}$, $t_1, t_2, t_3 \in \Z$ such that $G_{n,s} \simeq \langle \alpha, \beta \mid \alpha^2 = \beta^n = 1, \beta \alpha = \alpha \beta^s \rangle$, $\gcd(t_1-t_2,n) = 1$ and $S = (\beta^{t_1})^{[2n-1]} \bd (\beta^{t_2})^{[n-1]} \bd (\alpha \beta^{t_3})$.

\item Let $n$ be odd, and $|S| = \s(G_{n,s}) - 1 = \E(G_{n,s}) - 1 = 3n-1$. Suppose additionally that $n_1 > 3$, where $n_1$ is given by Lemma \ref{lemmafactorn}. Then $S$ is $2n$-product-one free if and only if  there exist $\alpha, \beta \in G_{n,s}$, $t_1, t_2, t_3 \in \Z$ such that $G_{n,s} \simeq \langle \alpha, \beta \mid \alpha^2 = \beta^n = 1, \beta \alpha = \alpha \beta^s \rangle$, $\gcd(t_1-t_2,n) = 1$ and $S = (\beta^{t_1})^{[2n-1]} \bd (\beta^{t_2})^{[n-1]} \bd (\alpha \beta^{t_3})$. 
\end{enumerate}
\end{theorem}

The assumption ``if $n$ is odd, then $n_1 > 3$'' of previous theorems relies on the fact that $D_6$ is isomorphic to some group quotient by Lemma \ref{lemmanormalsubgroups}{\it (a)}, therefore Lemma \ref{lemmaegzinversedihedral}{\it (c)} brings up some extra difficulties and the same arguments than all other cases do not work.

Fixed one of the invariants treated in this paper, the arguments for both direct and inverse results follow closely the same steps, in opposition to the direct and inverse problems related to $\d(G_{n,s})$. In fact, $\d(G_{n,s}) = n$ follows directly from a trivial upper bound, while its associated inverse problem (see \cite{MR3}) is pretty hard and follows from the study of many cases and subcases. For the problems dealt with here, it is also required to split into some cases depending on the parity of $n_1$ and $n_2$.


%

The paper is organized as follows. In Section \ref{lemmas}, we present some auxiliary results that will be used throughout the paper, which include the aforementioned factorization of $n$ depending on $s$, the inverse problems over cyclic and dihedral groups, some normal subgroups of $G_{n,s}$, and other generating sets. In Section \ref{sec:eta}, we prove Theorem \ref{thmeta} using the results on small Davenport constant. In Section \ref{sec:even}, we deal with Theorems \ref{thmesdirect} and \ref{thmesinverse} for $n$ even. In Section \ref{sec:odd}, we deal with the case $n$ odd.

\section{Auxiliary results}\label{lemmas}

In this section, we present some auxiliary results that will be used throughout the paper. The first one is a lemma that has been considered previously, and helps us to factorize $n$ nicely. But before we should notice that $s^2 \equiv 1 \pmod n$ and $s \not\equiv \pm1 \pmod n$ ensure that $n$ is neither an odd prime power nor twice an odd prime power.

\begin{lemma}[{\cite[Section~5]{GeGr}, \cite[Lemma 2.2]{MR3}}]\label{lemmafactorn}
Let $n \ge 8$ and $s$ be integers satisfying $s^2 \equiv 1 \pmod n$, but $s \not\equiv \pm1 \pmod n$.
\begin{enumerate}[(a)]
\item If both $n \neq p^t$ and $n \neq 2p^t$ for every prime $p$ and every integer $t \ge 1$, then there exist coprime integers $n_1, n_2 \ge 3$ such that $s \equiv -1 \pmod {n_1}$, $s \equiv 1 \pmod {n_2}$, and either 
$n = n_1n_2$ or 
$n = 2n_1n_2$.
\item If $n = 2^t$ for some $t \ge 3$, then 
$n = 2n_1n_2$, where either $(n_1,n_2) = (1,2^{t-1})$ satisfies $s \equiv 1 \pmod {n_2}$ or $(n_1,n_2) = (2^{t-1},1)$ satisfies $s \equiv -1 \pmod {n_1}$.
\end{enumerate}
\end{lemma}


The next results concern on the inverse problem related to Erd\H os-Ginzburg-Ziv constant over cyclic and dihedral groups, respectively.

\begin{lemma}[{\cite{Ga1}}]\label{lemmaegzinversecyclic}
Let $n \ge 2$ be an integer and let $2 \le k \le \lfloor n/2 \rfloor + 2$. Suppose that $S \in \F(C_n)$ and $|S| = 2n-k$. If $S$ is $n$-product-one free, then there exists $a \bd b \mid S$ such that $\min\{v_a(S), v_b(S)\} \ge n - 2k + 3$ and $v_a(S) + v_b(S) \ge 2n - 2k + 2$, where $ab^{-1}$ generates $C_n$.
\end{lemma}

\begin{lemma}[{\cite[Theorem 1.2]{OhZh}}]\label{lemmaegzinversedihedral}
Let $n \ge 3$ be an integer and consider $S \in \F(D_{2n})$. 
\begin{enumerate}[(a)]
\item If $n \ge 4$ is even and $|S| = \s(D_{2n}) - 1 = 2n - 1$, then $S$ is $n$-product-one free if and only if there exist $\alpha, \beta \in D_{2n}$ and $t_1, t_2, t_3 \in \Z$ such that $G_{n,s} \simeq \langle \alpha, \beta \mid \alpha^2 = \beta^n = 1, \beta \alpha = \alpha \beta^{-1} \rangle$, $\gcd(t_1-t_2,n) = 1$ and $S = (\beta^{t_1})^{[n-1]} \bd (\beta^{t_2})^{[n-1]} \bd (\alpha \beta^{t_3})$.
\item If $n \ge 4$ and $|S| = \E(D_{2n}) - 1 = 3n - 1$, then $S$ is $2n$-product-one free if and only if there exist $\alpha, \beta \in D_{2n}$ and $t_1, t_2, t_3 \in \Z$ such that $G_{n,s} \simeq \langle \alpha, \beta \mid \alpha^2 = \beta^n = 1, \beta \alpha = \alpha \beta^{-1} \rangle$, $\gcd(t_1-t_2,n) = 1$ and $S = (\beta^{t_1})^{[2n-1]} \bd (\beta^{t_2})^{[n-1]} \bd (\alpha \beta^{t_3})$.
\item If $n = 3$ and $|S| = \E(D_6) - 1 = 8$, then $S$ is $6$-product-one free if and only if there exist $\alpha, \beta \in D_6$ and $t_1, t_2, t_3 \in \Z$ such that $G_{n,s} \simeq \langle \alpha, \beta \mid \alpha^2 = \beta^3 = 1, \beta \alpha = \alpha \beta^{-1} \rangle$, $\gcd(t_1-t_2,3) = 1$ and either $S = (\beta^{t_1})^{[5]} \bd (\beta^{t_2})^{[2]} \bd (\alpha \beta^{t_3})$ or $S = 1^{[5]} \bd \alpha \bd \alpha \beta \bd \alpha \beta^2$.
\end{enumerate}
\end{lemma}

The following lemma provides the solution of the inverse problem related to small Davenport constant over $G_{n,s}$.

\begin{lemma}[{\cite[Theorem 1.1]{MR3}}]\label{lemmadavinversegns}
Let $S \in \F(G_{n,s})$ of length $\d(G_{n,s}) = n$. Then $S$ is product-one free if and only if there exist $\alpha, \beta \in G_{n,s}$ such that $G_{n,s} \simeq \langle \alpha, \beta \mid \alpha^2 = \beta^n = 1, \beta \alpha = \alpha \beta^s \rangle$ and $S = \beta^{[n-1]} \bd (\alpha \beta^t)$ for $0 \le t \le n-1$.
\end{lemma}

Sometimes, it is easier to obtain disjoint subproducts that belong to subgroups and then it is required to know some normal subgroups of $G_{n,s}$. The following lemma provides what is needed.

\begin{lemma}\label{lemmanormalsubgroups}
Let $x, y \in G_{n,s}$ be as in Equation \eqref{Gns}, where either $n = n_1n_2$ or $n = 2n_1n_2$ are as in Lemma \ref{lemmafactorn}.
\begin{enumerate}[(a)]
\item Suppose that $n = hn_1n_2$ for $h \in \{1,2\}$ and $n_2$ is odd. Then $H = \langle y^{n_1} \rangle \simeq C_{hn_2}$ is a normal subgroup of $G_{n,s}$ and satisfies $G_{n,s}/H \simeq D_{2n_1} \simeq \langle x,y^{hn_2} \rangle$.

\item Suppose that $n = hn_1n_2$ for $h \in \{1,2\}$ and $n_2$ is even. Then $H = \langle y^{2n_1} \rangle \simeq C_{hn_2/2}$ is a normal subgroup of $G_{n,s}$ and satisfies $G_{n,s}/H \simeq D_{4n_1} \simeq \langle x,y^{hn_2} \rangle$. \\
In the case that $n_1 = 1$, $D_{4n_1} \simeq C_2^2$ is the Klein group.

\item Suppose that $n = n_1n_2$. Then $H = \langle x,y^{n_2} \rangle \simeq D_{2n_1}$ is a normal subgroup of $G_{n,s}$ and satisfies $G_{n,s}/H \simeq C_{n_2} \simeq \langle y^{n_1} \rangle$. 
\end{enumerate}
\end{lemma}

\proof
We will prove only {\it (a)} for $n = n_1n_2$; the other cases follow similarly. Let $x^{\alpha}y^{\beta} \in G_{n,s}$ and $y^{tn_1} \in H$. We have 
$$x^{\alpha}y^{\beta} \cdot y^{tn_1} \cdot (x^{\alpha}y^{\beta})^{-1} = y^{tn_1} \in H$$
since $s \equiv 1 \pmod {n_2}$, 
thus 
$H$ is a normal subgroup of $G_{n,s}$. Furthermore, $x$ has order $2$, $y^{n_2}$ has order $n_1$, and $y^{n_2}x = xy^{n_2s} = xy^{-n_2} = x(y^{n_2})^{-1}$ since $s \equiv -1 \pmod {n_1}$. Hence, $G_{n,s}/H \simeq D_{2n_1}$.
\qed

\vspace{3mm}

Finally, the next lemma provides some generating sets other than $G_{n,s} = \langle x,y \rangle$ that we may consider.

\begin{lemma}\label{lemmagenerators}
Let $x, y \in G_{n,s}$ be as in Equation \eqref{Gns}. It holds $G_{n,s} \simeq 
\langle x, y^v \mid  \gcd(v,n) = 1 \rangle$. Furthermore, for $t \ge 3$, it holds $G_{2^t,2^{t-1}+1} \simeq \langle x, xy^v \mid v \text{ is odd} \rangle$.
\end{lemma}

\proof
The first part, $G_{n,s} \simeq \langle x, y^v \mid \gcd(v,n) = 1 \rangle$ follows from the fact that $y^v$ has order $n$ since $(y^v)^r = y^{rv} = 1$ if and only if $r \equiv 0 \pmod n$, and $y^v \cdot x = x \cdot (y^v)^s$.

Moreover, for $t \ge 3$ we have $G_{2^t,2^{t-1}+1} \simeq \langle x, xy^v \mid v \text{ is odd} \rangle$. Indeed, $xy^v$ has order $2^t$ since $(xy^v)^r = 1$ implies that $r$ is even, say $r = 2r'$, thus $(xy^v)^r = (y^{2v + 2^{t-1}})^{r'} = (y^{v+2^{t-2}})^r = 1$ if and only if $r \equiv 0 \pmod {2^t}$, and then $xy^v \cdot x = y^{v(2^{t-1}+1)} = x \cdot (xy^v)^{2^{t-1}+1}$.
\qed

\section{On the $\eta$-constant: Proof of Theorem \ref{thmeta}}\label{sec:eta}

Let $H = \langle y \rangle$, where $x,y \in G_{n,s}$ are as in Equation \eqref{Gns}. For $n$ odd, we have $\exp(G_{n,s}) = 2n > n = \d(G_{n,s})$, therefore $\e(G_{n,s}) = n+1$. Now we let $n$ be even, thus $\exp(G_{n,s}) = n$. In order to show that $\e(G_{n,s}) = n+1$, it is required to show that any sequence of length $n+1$ has some short product-one subsequence. Notice that the sequences given by Lemma \ref{lemmadavinversegns} have only one term out of $H$ and $n-1$ terms in $H$. Hence let $S \in \F(G_{n,s})$ with $|S| = n+1$. We will remove a term $g \mid S$ in such way that $S \bd g^{[-1]}$ (that has length $n$) is not product-one free. If $1 \le |S \cap H| \le n+1$, then there exists $g \mid S$ such that $|(S \bd g^{[-1]}) \cap H|$ is even, thus it can not be one. By Lemma \ref{lemmadavinversegns}, $S \bd g^{[-1]}$ is not product-one free, therefore $S$ has a short product-one subsequence. If $S \cap H$ is the empty sequence, then one can remove any term $g \mid S$, so that $|(S \bd g^{[-1]}) \cap H| = 0$ and $S \bd g^{[-1]}$ has length $n$, hence $S \bd g^{[-1]}$ is product-one free if and only if $S \bd g^{[-1]} = (xy^v)^{[n-1]} \bd (xy^u)$, where $v$ is odd and $u$ is even, by Lemmas \ref{lemmadavinversegns} and \ref{lemmagenerators}. In this case, one could remove another term by considering that either $S \bd (xy^v)^{[-1]}$ or $S \bd (xy^u)^{[-1]}$ is not of the later form. Thus $S$ has a short product-one subsequence.

On the other hand, the sequences of length $n$ which are product-one free over $G_{n,s}$ are precisely the sequences of length $n$ which have no short product-one subsequences, and are given by Lemma \ref{lemmadavinversegns}.
\qed

\section{On Erd\H os-Ginzburg-Ziv direct and inverse theorems: The case $n$ even}

Let $n$ be even. The lower bounds $\s(G_{n,s}) \ge 2n$ and $\E(G_{n,s}) \ge 3n$ follow from inequalities \eqref{lowerboundse}, therefore it is only required to obtain the respective upper bounds. In view of Lemma \ref{lemmafactorn}, we split the proof into the following cases:
\begin{enumerate}[(i)]
\item $n=2n_1n_2$ (similarly $n=n_1n_2$), with $n_1,n_2\geq3$ coprime integers and exactly one being even; 
\item $n=2^t=2n_1n_2$, with $t \geq3$, and $(n_1,n_2)=(2^{t-1},1)$; or
\item $n=2^t=2n_1n_2$, with $t\geq3$, and $(n_1,n_2)=(1,2^{t-1})$.
\end{enumerate}

Moreover, the direct and inverse problems in this section are solved in a similar way.

\subsection{The direct problems: Proof of Theorem \ref{thmesdirect} for $n$ even}\label{sec:even}

Let $S \in \F(G_{n,s})$ with $|S|=2n$. First of all, we are going to show that $\s(G_{n,s}) = 2n$.

The proof for cases (i) with $n_1$ even and (ii) follow the same steps. In fact, assume that $n_1$ is even. 
Let $n=2n_1n_2$ (resp. $n=n_1n_2$) and let $H=\langle y^{n_1}\rangle\simeq C_{2n_2}$ (resp. $H\simeq C_{n_2}$) be a normal subgroup of $G_{n,s}$, so that $\overline{G_{n,s}} = G_{n,s}/H \simeq D_{2n_1}$ by Lemma \ref{lemmanormalsubgroups}. 
Since $\s(\overline{G_{n,s}})=2n_1$ and $|\overline{S}|=2n>2n_1$, $\overline{S}$ contains some $n_1$-product-one subsequence $\overline{T}_1$, that is, denoting by $T_1$ the subsequence of $S$ formed by the pre-image by $\phi: G_{n,s} \to \overline{G_{n,s}}$ of the terms of $\overline{T}_i$, it follows that $\pi(T_1) \subset H$ since $H$ is a normal subgroup. Since $|S \bd T_1^{[-1]}| = 2n - n_1$, we construct inductively disjoint subsequences $T_2, \dots, T_{\ell}$ of $S$ with $|T_i| = n_1$ and $\pi(T_i) \in H$ until $$|\overline{S}{\bd}\overline{T}_1^{[-1]}{\bd}\cdots{\bd}\overline{T}_{\ell}^{[-1]}| = 2n-n_1{\ell}<2n_1,$$ that is, $4n_1n_2-n_1{\ell}<2n_1$ (resp. $2n_1n_2-n_1{\ell}<2n_1$), which implies that ${\ell}=4n_2-1$ (resp. ${\ell}=2n_2-1$). 
Since $\s(H)=4n_2-1$ (resp. $\s(H)=2n_2-1$), there are $j=2n_2$ (resp. $j=n_2$) subsequences $\overline{T}_{i_1},\dots,\overline{T}_{i_j}$ such that $1 \in \pi(T_{i_1} \bd \dots \bd T_{i_j})$. Anyway, $T_{i_1}{\bd}\cdots{\bd}T_{i_j} \mid S$ is a $n$-product-one subsequence, thus $\s(G_{n,s})=2n$.

We now assume $n_2$ even in case (i). Let $n=2n_1n_2$ (resp. $n=n_1n_2$) and let $H = \langle y^{2n_1}\rangle\simeq C_{n_2}$ (resp. $H \simeq C_{n_2/2}$) and $\overline{G_{n,s}}=G_{n,s}/H\simeq D_{4n_1}$. We emphasize that the only differences from the previous arguments for this one is that $\s(\overline{G_{n,s}})=4n_1$ and $\s(H)=2n_2-1$ (resp. $\s(H)=n_2-1$). Therefore, following the same arguments we prove that $\s(G_{n,s})=2n$.

We now assume case (iii). In particular, $s=2^{t-1}+1$. Suppose that $S \in \F(G_{n,s})$ is $n$-product-one free, and fix $H = \langle y^{2}\rangle\simeq C_{n/2}=C_{2^{t-1}}$ so that $\overline{G_{n,s}}=G_{n,s}/H\simeq C_2\times C_2$. 
Since $\s(\overline{G_{n,s}})=5$ and $|\overline{S}|=2^{t+1}>5$, $\overline{S}$ contains a $2$-product-one subsequence $\overline{T}_1$. This argument can be repeated to construct disjoint $2$-product-one subsequences $\overline{T}_2, \dots, \overline{T}_{\ell}$ of $\overline{S} \in \F(\overline{G_{n,s}})$ with $1 \in \pi(\overline{T}_i)$ until 
$$|\overline{S} {\bd}\overline{T}_1^{[-1]}{\bd}\cdots{\bd}\overline{T}_{\ell}^{[-1]}| = 2^{t+1}-2{\ell} < 5,$$ 
that is, ${\ell}=2^t-2$. 
Then $$S = T_1 \bd \dots \bd T_{2^t-2}\bd T,$$ where $T_i$ are subsequences of length $2$ such that $\pi(T_i) \subset H$ and, since $S$ is $n$-product-one free, $\overline{T}$ is also $2$-product-one free with $|T|=4$. It implies that $T=y^{2\alpha}\bd y^{2\beta+1}\bd y^{2\gamma}x\bd y^{2\delta+1}x$.
Let $h_i \in \pi(T_i)$. Since $\s(H)=2^t-1$, it follows from Lemma \ref{lemmaegzinversecyclic} that $$h_1\bd\cdots\bd h_{2^t-2}=(y^{2a})^{[2^{t-1}-1]}\bd(y^{2b})^{[2^{t-1}-1]},$$ where $a-b$ is odd. 
We claim that $S$ is not $2^t$-product-one free. In fact, let $k \in [0,2^{t-1}-1]$ such that 
$$y^{2\alpha} \cdot y^{2\beta+1} \cdot (y^{2a})^k \cdot y^{2\gamma}x \cdot (y^{2b})^{2^{t-1}-k-2} \cdot y^{2\delta+1}x = 1.$$ 
Such $k$ exists since the previous equation is equivalent to 
$$2\alpha+2\beta+1+2ak+2\gamma-2bks-4bs+2\delta s+s\equiv 0 \!\!\! \pmod {2^t}.$$
Since $s=2^{t-1}+1$, the latter holds if and only if 
$$(a-b)k\equiv 2b-1-\alpha-\beta-\gamma-\delta-2^{t-2} \!\!\! \pmod {2^{t-1}}.$$ 
Since $\gcd(a-b,2^{t-1})=1$, it has a solution $k \in [0, 2^{t-1}-1]$. If $k\in[0,2^{t-1}-2]$, then $S$ is not $2^t$-product-one free and we are done. On the other hand, if $k=2^{t-1}-1$, then we have
\begin{equation}\label{eq:u+v}
a+b \equiv 2^{t-2}+1+\alpha+\beta+\gamma+\delta \!\!\! \pmod {2^{t-1}}.
\end{equation}
In this case, there exists $m \in [0,2^{t-1}-2]$ such that 
$$y^{2\alpha} \cdot y^{2\beta+1} \cdot (y^{2a})^{m} \cdot y^{2\delta+1}x \cdot (y^{2b})^{2^{t-1}-m-2} \cdot y^{2\gamma}x = 1.$$ 
In fact, the previous equation holds if and only if
$$(a-b) m \equiv 2b-1-\alpha-\beta-\gamma-\delta \!\!\! \pmod {2^{t-1}},$$
which has a solution $m \in [0,2^{t-1}-1]$. Again we are done if $m \in [0,2^{t-1}-2]$. If $m = 2^{t-1}-1$, then $a+b \equiv 1+\alpha+\beta+\gamma+\delta \pmod {2^{t-1}}$, which contradicts Eq. \ref{eq:u+v}. Therefore, $\s(G_{n,s})=2n$.

Now, we conclude that $\E(G_{n,s}) \le 3n$ by extracting from any $S \in \F(G_{n,s})$ with $|S| = 3n$ two disjoint $n$-product-one subsequences. Therefore, $\E(G_{n,s}) = 3n$.
\qed

\subsection{The inverse problems: Proof of Theorem \ref{thmesinverse}(a)}

In order to solve the inverse problem, let us first assume that $S$ is $n$-product-one free with $|S|=2n-1$, $n_1$ is even and either (i) or (ii) holds. Similar to the direct problem for $n=2n_1n_2$ (resp. $n=n_1n_2$), we obtain that $\overline{S}=\overline{T}_1{\bd}\cdots{\bd}\overline{T}_{4n_2-2}\bd\overline{T}$ (resp. $\overline{S}=\overline{T}_1{\bd}\cdots{\bd}\overline{T}_{2n_2-2}\bd\overline{T}$), with $|T_i| = n_1$, $\pi(T_i) \subset H \simeq C_{2n_2}$ (resp. $H \simeq C_{n_2}$), and $|\overline{T}|=2n_1-1$. 
If $\overline{T}$ contains a $n_1$-product-one subsequence, then $S$ contains a $n$-product-one subsequence. Hence, $\overline{T}$ must be $n_1$-product-one free. It follows from Lemma \ref{lemmaegzinversedihedral} that 
$$\overline{T} = (\overline{y^{a}})^{[n_1-1]} \bd (\overline{y^{b}})^{[n_1-1]}\bd\overline{y^{r_1} x},$$ 
where $\gcd(a-b,n_1)=1$. For each $i \in [1,4n_2-2]$ (resp. $i \in [1,2n_2-2]$) and each $a_{i_j} \mid \overline{T}_i$, we consider $\overline{A}_{i_j} = \overline{T}_i{\bd}(a_{i_j})^{[-1]}{\bd}\overline{T}$. Since $n_1 \geq 4$ and $|\overline{A}_{i_j}| = 3n_1-2 \geq 2n_1 = \s(\overline{G_{n,s}})$, then 
$$\overline{S}=(\overline{y^{a}})^{[u]} \bd (\overline{y^{b}})^{[v]}\bd\ldotprod_{1\leq i\leq \ell}\overline{y^{r_i} x},$$ 
where $u,v\geq n_1-1$ and $\ell \geq 1$. We will show now that $u=v=n-1$ and $\ell=1$.  Let $\overline{B}=\overline{S}{\bd}(\overline{y^{r_1} x})^{[-1]}$. Similarly as done before, $\overline{B}=\overline{B}_1{\bd}\cdots{\bd}\overline{B}_{4n_2-2}{\bd}\overline{C}$ (resp. $\overline{B}=\overline{B}_1{\bd}\cdots{\bd}\overline{B}_{2n_2-2}{\bd}\overline{C}$), where each $\overline B_i$ is a $n_1$-product-one subsequence, and $|\overline{C}|=2n_1-2$. By considering $\overline{C}{\bd}(\overline{y}^{r_1}\overline{x})$, since $\overline{S}$ is $n_2$-product-one free (resp. $2n_2$-product-one free), we obtain that every term of $\overline{C}$ belongs to $\{\overline{y^{a}},\overline{y^{b}}\}$. Now, for each $i \in [1,4n_2-2]$ (resp. $i \in [1,2n_2-2]$) and each $b_{i_j} \mid \overline{B}_i$, we consider $\overline{C}_{i_j} = \overline{B}_i {\bd} (b_{i_j})^{[-1]} {\bd} \overline{C}$. Again, $n_1 \geq 4$ and $|\overline{C}_{i_j}| = 3n_1-3 \geq 2n_1 = \s(\overline{G_{n,s}})$, therefore 
$$\overline{S}=(\overline{y^{a}})^{[u]}\bd(\overline{y^{b}})^{[v]}\bd\overline{y^{r_1} x},$$ 
where $u,v \geq n_1-1$. Thus $|S\cap\langle y\rangle|=2n-2$ and, since $\s(\langle y \rangle) = 2n-1$, we conclude that $$S=(y^{a})^{[n-1]} \bd (y^{b})^{[n-1]} \bd y^{r}x,$$ 
where $\gcd(a-b,n)=1$. It is easy to show that this sequence is $n$-product-one free.

In the case (i) for $n_2$ even, the proof is completely similar.

We now assume case (iii) and $|S|=2n-1$ a $n$-product-one free sequence, where $n = 2^t$. Similarly as in previous cases, $S = T_1\bd\dots\bd T_{2^t-2}\bd T$, with $|T_i| = 2$, $\pi(T_i) \subset H = \langle y^2 \rangle$, and $|T|=3$. Let $h_i \in \pi(T_i)$. Again, it follows from Lemma \ref{lemmaegzinversedihedral} that 
\begin{equation}\label{h1hn}
h_1\bd\cdots\bd h_{2^t-2}=(y^{2a})^{[2^{t-1}-1]}\bd(y^{2b})^{[2^{t-1}-1]}, 
\end{equation}
where $a-b$ is odd (say, $a$ is even and $b$ is odd), and 
$$T=(y^{2\alpha})^{[\varepsilon]}\bd(y^{2\beta+1})^{[\lambda]}\bd(y^{2\gamma}x)^{[\theta]}\bd(y^{2\delta+1}x)^{[\nu]},$$ 
where $\varepsilon, \lambda, \theta, \nu \in \{0,1\}$ and $\varepsilon +  \lambda + \theta + \nu = 3$. For $i \in [1,2^t-2]$, let $T_i = g_{i,1} \bd g_{i,2}$. Since $\overline{g}_{i,1}\overline{g}_{i,2} = \overline 1$ in $G/H \simeq C_2 \times C_2$, if $g_{i,1}= y^{a_1}x^{b_1}$ and $g_{i,2}= y^{a_2}x^{b_2}$, then $a_1\equiv a_2 \pmod 2$ and $b_1\equiv b_2\pmod 2$.
Therefore, 
$$S = \ldotprod_{1\le i \le m_1} y^{2\alpha_i} \bd \ldotprod_{1\le i \le m_2} y^{2\beta_i+1} \bd \ldotprod_{1\le i \le m_3} y^{2\gamma_i}x \bd \ldotprod_{1\le i \le m_4} y^{2\delta_i+1}x,$$ 
with $m_1+m_2+m_3+m_4=2^{t+1}-1$ exactly one of the $m_i$'s being even (the product is empty if $m_i=0$). If $m_i \geq 1$ for every $i \in \{1,2,3,4\}$, let $y^{2k} = y^{2\gamma_1}x \cdot y^{2\delta_1+1}x \cdot y^{2\alpha_1} \cdot y^{2\beta_1+1}$ and $y^{2k'} = y^{2\delta_1+1}x \cdot y^{2\gamma_1}x \cdot y^{2\alpha_1} \cdot y^{2\beta_1+1} = y^{2k+2^{t-1}}$. Then either $y^{2k} \cdot (y^{2a})^r \cdot (y^{2b})^{2^{t-1}-r-2} = 1$ for some $r \in [0,2^{t-1}-2]$, or $y^{2k'} \cdot (y^{2a})^r \cdot (y^{2b})^{2^{t-1}-r-2} = 1$ for some $r \in [0,2^{t-1}-2]$, since $2k'\equiv 2k+2^{t-1} \pmod {2^t}$ and $\gcd(a-b,2^{t-1})=1$. Hence, we may assume that $m_i = 0$ for exactly one $i \in \{1,2,3,4\}$, thus there are four possibilities for $S$.


\begin{enumerate}[(1)]
\item $m_4 = 0$, that is, $S=\displaystyle\ldotprod_{1\le i \le 2m_1+1} y^{2\alpha_i} \bd \ldotprod_{1\le i \le 2m_2+1} y^{2\beta_i+1} \bd \ldotprod_{1\le i \le 2m_3+1} y^{2\gamma_i}x$. \\
We rewrite $S$ as 
$$\ldotprod_{1\le i \le m_{1,1}}y^{4\alpha_{i,0}}\bd\!\!\!\ldotprod_{1\le i \le m_{1,2}}y^{4\alpha_{i,2}+2}\bd\!\!\!\ldotprod_{1\le i \le m_{2,1}}y^{4\beta_{i,1}+1}\bd\!\!\!\ldotprod_{1\le i \le m_{2,2}}y^{4\beta_{i,3}+3}\bd\!\!\!\ldotprod_{1\le i \le m_{3,1}}y^{4\gamma_{i,0}}x\bd\!\!\!\ldotprod_{1\le i \le m_{3,2}}y^{4\gamma_{i,2}+2}x.$$  
Suppose that $y^{2a} = y^{4\beta_{i,1}+1} \cdot y^{4b_{j,3}+3}$. We are going to replace this term by another one, which contradicts the form given by \eqref{h1hn}, hence $S$ will not be $n$-product-one free. Since $T=y^{2a}\bd y^{2b+1}\bd y^{2c}x$, $u$ is even and $v$ is odd, either $y^{4\beta_{i,1}+1} \cdot y^{2b+1}$ or $y^{4\beta_{i,3}+3} \cdot y^{2b+1}$ belongs to $\langle y^2 \rangle \setminus \langle y^4 \rangle$, therefore it can not be $y^{2a} \in \langle y^4 \rangle$. This argument can be applied to show that $\alpha_{i,0} \equiv a_{j,0} \pmod {n/4}$, $\alpha_{i,2} \equiv a_{j,2} \pmod {n/4}$, $\beta_{i,1} \equiv b_{j,1} \pmod {n/4}$, $\beta_{i,3} \equiv b_{j,3} \pmod {n/4}$, $\gamma_{i,0} \equiv \gamma_{j,0} \pmod {n/4}$, and $\gamma_{i,2} \equiv \gamma_{j,2} \pmod {n/4}$ for every $i,j$. Using the same argument for $y^{2b}$, we obtain 
\begin{align*}
y^{2a} &\in \{(y^{4\alpha_{i,0}})^2, (y^{4\alpha_{i,2}+2})^2, (y^{4\gamma_{i,0}}x)^2, (y^{4\gamma_{i,2}+2}x)^2\}, \\
y^{2b} &\in \{(y^{4\beta_{i,1}+1})^2, (y^{4\beta_{i,3}+3})^2\}.
\end{align*} 
If there exist $i,j$ such that $y^{2b} = (y^{4\beta_{i,1}+1})^2 = (y^{4b_{j,3}+3})^2$, then we replace these products by two of the form $y^{4\beta_{i,1}+1} \cdot y^{4b_{j,3}+3}$, both belonging to $\langle y^4 \rangle$, which contradicts the form given by \eqref{h1hn}. Thus $S$ is not $n$-product-one free. It implies that $2m_2 + 1 = 2^t - 1$ and either $m_{2,1}=0$ or $m_{2,2}=0$. Similarly, either $m_{1,1}=0$ or $m_{1,2}=0$, and either $m_{3,1}=0$ or $m_{3,2}=0$. Furthermore, $2(m_1 + m_3) + 1 = 2^t - 1$. \\
Assume $m_3 \geq 1$. Let $y^{2k} = y^{2b_1+1} \cdot y^{2c_1}x \cdot y^{2c_2}x \cdot y^{2b_2+1}$ and $y^{2k'} = y^{2c_1}x \cdot y^{2b_1+1} \cdot y^{2c_2}x \cdot y^{2b_2+1} = y^{2k+2^{t-1}}$. Using the same argument as in Subsection \ref{sec:even}(iii), we obtain a contradiction. Therefore $m_3=0$, and $$S=(y^{2a})^{[2^t-1]}\bd(y^{2b+1})^{[2^t-1]}\bd y^{2c}x.$$
Furthermore, it is easy to check that the latter is $n$-product-one free.

\item $m_3 = 0$, that is, $\displaystyle S=\ldotprod_{1\le i \le 2m_1+1} y^{2\alpha_i}\bd\ldotprod_{1\le i \le 2m_2+1} y^{2\beta_i+1}\bd\ldotprod_{1\le i \le 2m_4+1} y^{2\delta_i+1}x$. \\
We rewrite $S$ as 
$$\ldotprod_{1\le i \le m_{1,1}}y^{4\alpha_{i,0}}\bd\!\!\!\ldotprod_{1\le i \le m_{1,2}}y^{4\alpha_{i,2}+2}\bd\!\!\!\ldotprod_{1\le i \le m_{2,1}}y^{4\beta_{i,1}+1}\bd\!\!\!\ldotprod_{1\le i \le m_{2,2}}y^{4\beta_{i,3}+3}\bd\!\!\!\ldotprod_{1\le i \le m_{3,1}}y^{4\delta_{i,1}+1}x\bd\!\!\!\ldotprod_{1\le i \le m_{3,2}}y^{4\delta_{i,3}+3}x.$$ 
By the same arguments as previous case, we obtain $\alpha_{i,0} \equiv a_{j,0} \pmod {n/4}$, $\alpha_{i,2} \equiv a_{j,2} \pmod {n/4}$, $\beta_{i,1} \equiv b_{j,1} \pmod {n/4}$, $\beta_{i,3} \equiv b_{j,3} \pmod {n/4}$, $\delta_{i,1} \equiv \delta_{j,1} \pmod {n/4}$, and $\delta_{i,3} \equiv \delta_{j,3} \pmod {n/4}$ for every $i,j$. Furthermore,
\begin{align*}
y^{2a} &\in \{(y^{4\alpha_{i,0}})^2, (y^{4\alpha_{i,2}+2})^2\} \\
y^{2b} &\in \{(y^{4\beta_{i,1}+1})^2, (y^{4\beta_{i,3}+3})^2, (y^{4\delta_{i,1}+1}x)^2, (y^{4\delta_{i,3}+3}x)^2\}.
\end{align*}
Moreover, either $m_{1,1}=0$ or $m_{1,2}=0$ (which implies that $2m_1+1=2^t-1$), either $m_{2,1}=0$ or $m_{2,2}=0$, and either $m_{4,1}=0$ or $m_{4,2}=0$. If $m_2 \ge 1$ and $m_4 \geq 1$, we set $y^{2k}=y^{2b_1+1} \cdot y^{2d_1+1}x \cdot y^{2b_2+1} \cdot y^{2d_2+1}x$ and $y^{2k'} = y^{2b_1+1} \cdot y^{2b_2+1} \cdot y^{2d_1+1}x \cdot y^{2d_2+1}x = y^{2k+2^{t-1}}$. The same argument as Subcase (1) implies that $S$ is not $n$-product-one free, therefore either $m_2=0$ and $2m_4+1=2^t-1$, or $m_4 = 0$ and $2m_2+1=2^t-1$. It implies that either 
$$S=(y^{2a})^{2^t-1}\bd(y^{2b+1})^{2^t-1}\bd y^{2d+1}x \quad \text{ or } \quad S=(y^{2a})^{2^t-1}\bd y^{2b+1}\bd(y^{2d+1}x)^{2^t-1}.$$
Furthermore, it is easy to check that the latter sequences are $n$-product-one free.


\item $m_2 = 0$, that is, $\displaystyle S=\ldotprod_{1\le i \le 2m_1+1} y^{2\alpha_i} \bd \ldotprod_{1\le i \le 2m_3+1} y^{2\gamma_i}x \bd \ldotprod_{1\le i \le 2m_4+1} y^{2\delta_i+1}x$. \\
This case is completely similar and returns the $n$-product-one free sequences $$S=(y^{2a})^{2^t-1}\bd y^{2c}x\bd(y^{2d+1}x)^{2^t-1}.$$


\item $m_1 = 0$, that is, $\displaystyle S=\ldotprod_{1\le i \le 2m_2+1} y^{2\beta_i+1}\bd\ldotprod_{1\le i \le 2m_3+1} y^{2\gamma_i}x\bd\ldotprod_{1\le i \le 2m_4+1} y^{2\delta_i+1}x.$ \\
This case is completely similar and returns the $n$-product-one free sequences 
$$S = y^{2b+1} \bd (y^{2c}x)^{[2^{t+1}-3]} \bd y^{2d+1}x.$$ 
However, this $S$ is not $n$-product-one free, since $(y^{2c}x)^{2^t}=1$.

\end{enumerate}
\qed

\subsection{The inverse problems: Proof of Theorem \ref{thmesinverse}(b)}

Let $S \in \F(G_{n,s})$ be a $2n$-product-one free sequence of length $|S| = 3n-1$. Since $|S| > 2n = \s(G_{n,s})$, it is possible to extract $n$-product-one subsequence $T \mid S$. It follows that $|S \bd T^{[-1]}| = 2n-1$ and $S \bd T^{[-1]}$ is $n$-product-one free. By Theorem \ref{thmesinverse}(a), it follows that $S \bd T^{[-1]} = (y^{t_1})^{[n-1]} \bd (y^{t_2})^{[n-1]} \bd xy^{t_3}$, where $\gcd(t_1-t_2,n) = 1$.

Suppose that $S$ contains two terms, say $g_1$ and $g_2$, which are both distinct from $y^{t_1}$ and $y^{t_2}$. Then $|S \bd (g_1 \bd g_2)^{[-1]}| = 3n-3 > 2n$, therefore we may apply the same argument as before to ensure that there exists a $n$-product-one subsequence $T \mid (S \bd (g_1 \bd g_2)^{[-1]}$, however the remaining $|S \bd T^{[-1]}| = 2n-1$ terms have at least two terms other than $y^{t_1} \bd y^{t_2}$. By Lemma \ref{lemmaegzinversedihedral}, $S$ contains a $2n$-product-one subsequence in this case.

It implies that $S = (y^{t_1})^{[u]} \bd (y^{t_2})^{[3n-2-u]} \bd xy^{t_3}$. Since $\gcd(t_1-t_2,n) = 1$, it follows that the only ways of obtaining a product-one subsequence of length $n$ are $(y^{t_1})^n = 1$ or $(y^{t_2})^n = 1$. Therefore, either $u = 2n-1$ or $u = n-1$, and we are done.
\qed

\section{On Erd\H os-Ginzburg-Ziv direct and inverse theorems: The case $n$ odd}\label{sec:odd}

Throughout this section, we assume that $n = n_1n_2$ is odd. We first consider the direct problem (Theorem \ref{thmesdirect}) in Subsection \ref{subsectiondirectnodd} and then the inverse problem (Theorem \ref{thmesinverse}(c)) in Subsection \ref{subsectioninversenodd}. It is worth mentioning that $\exp(G_{n,s}) = |G_{n,s}| = 2n$ in this case, thus $\s(G_{n,s}) = \E(G_{n,s})$. The missing case is $n_1 = 3$, for which there are some extra difficulties in view of Lemmas \ref{lemmaegzinversedihedral}(c) and \ref{lemmanormalsubgroups}, thus our arguments do not work. It will be useful to shorten our notation by setting $S_y = S \cap \langle y\rangle$, $S_{y^{n_1}} = S \cap \langle y^{n_1}\rangle$, and $S_{xy} = S \cap (\langle y\rangle x)$.

\subsection{The direct problem: Proof of Theorem \ref{thmesdirect} for $n$ odd and $n_1>3$}\label{subsectiondirectnodd}

The lower bound $\s(G_{n,s}) \ge 3n$ follows from the first inequality of \eqref{lowerboundse}, therefore it is just required to prove the upper bound. 
Let $n$ be odd and assume that $S \in \F(G_{n,s})$ has length $|S|=3n$ and is $2n$-product-one free. 

Since $\s(\langle y\rangle) = 2n-1$, if $|S_y| \ge 3n-1$, then $S$ is not $2n$-product-one free because it is possible to select two disjoint $n$-product-one subsequences. Therefore, we will assume that $|S_{xy}| \geq 2$.

Let $H=\langle x,y^{n_2}\rangle\simeq D_{2n_1}$, so that $\overline{G_{n,s}} = G_{n,s}/H \simeq C_{n_2}$. 
Since $\s(\overline{G_{n,s}})=2n_2-1$ and $|\overline{S}|=3n>2n_2-1$, $\overline{S}$ contains a $n_2$-product-one subsequence $\overline{T}_1$. In other words, there exists $T_1 \mid S$ with $|T_1| = n_2$ such that $\pi(T_1) \subset H$. This argument can be applied to construct similarly disjoint subsequences $T_2, \dots, T_{3n_1-1}$ of $S$, each of length $n_2$, whose products belong to $H$. For each $i \in [1, 3n_1-1]$, let
$${T}_i={g}_{i,1}{\bd}\cdots{\bd}{g}_{i,n_2} \quad \text{ and } \quad {h}_i={g}_{i,1}\cdots{g}_{i,n_2}\in \pi(T_i) \subset H.$$

Since $\s(H)=3n _1$, $n_1>3$ and $S$ is $2n$-product-one free, it follows from Lemma \ref{lemmaegzinversedihedral} that 
\begin{equation}\label{sequencehi}
{h}_1\bd\cdots\bd{h}_{3n_1-1}=(y^{an_2})^{[2n_1-1]}\bd ({y}^{bn_2})^{[n_1-1]}\bd {y}^{r_1n_2}{x} \in \F(H),
\end{equation} 
where $\gcd (a-b,n_1)=1$.

Let $S = T_1\bd\cdots\bd T_{3n_1-1}\bd\ldotprod_{1\leq i\leq n_2}e_i$, where $\pi(e_1 \bd \dots \bd e_{n_2})$ does not contain any element of $H$. Reordering the indexes, we can write 
$$\begin{array}{ll}h_i=y^{an_2}&\mbox{ if }i\in[1,2n_1-1],\\h_i=y^{bn_2}&\mbox{ if }i\in[2n_1,3n_1-2],\\h_i=y^{r_1n_2}x&\mbox{ if }i=3n_1-1.
\end{array}$$

In particular, for $i \in [1,3n_1-2]$, the $h_i$'s must be unique, that is, $\pi(T_i) = \{h_i\}$; in other words, the terms of $T_i$ must commute between them. Notice that the product of two terms in $\langle y \rangle$ always commute. On the other hand, $y^{\alpha} \cdot y^{\beta} x = y^{\beta}x \cdot y^{\alpha}$ if and only if $\alpha \equiv 0 \pmod {n_1}$ since $\gcd(n_1,n_2)=1$, $s \equiv -1 \pmod {n_1}$, $s \equiv 1 \pmod {n_2}$, and $n_1$ is odd. Furthermore, $y^{\alpha}x \cdot y^{\beta}x = y^{\beta}x \cdot y^{\alpha}x$ if and only if $\alpha \equiv \beta \pmod {n_1}$ by the same reason. 

Moreover, if $i \in [1,3n_1-2]$ and $g_{i,j} \in \langle y \rangle x$ for some $j \in [1,n_2]$, then $h_i = 1$. Indeed, since they commute, we may assume $h_i = (y^{\beta_{i,1}}x \dots y^{\beta_{i,2k}}x) \cdot (y^{\alpha_{i,2k+1}} \dots y^{\alpha_{i,n_2}}) = y^{tn_2}$ (the order of $x$ forces an even number of $j \in [1,n_2]$ with $g_{i,j} \in \langle y \rangle x$). Looking at the exponent of $y$ in $h_i$, it follows that $\beta_{i,1} + \dots + \beta_{i,2k} + \alpha_{i,2k+1} + \dots + \alpha_{i,n_2} \equiv tn_2 \equiv 0 \pmod {n_2}$ since $s \equiv 1 \pmod {n_2}$, and $tn_2 \equiv (\beta_{i,1} - \beta_{i,2}) + \dots + (\beta_{i,2k-1} - \beta_{i,2k}) +  \alpha_{i,2k+1} + \dots + \\alpha_{i,n_2} \equiv 0 \pmod {n_1}$ since $s \equiv -1 \pmod {n_1}$, thus $t \equiv 0 \pmod {n_1}$. Therefore, $h_i = 1$.

We now split into some cases depending on $|S_{xy}|$.


\begin{enumerate}[(1)]
\item Case $|S_{xy}| \geq 3n_2-1$. In this case, it is possible to obtain two $n_2$-product-one subsequences of $\overline{S_{xy}}$. Since $n_2$ is odd, each of these products belongs to $\langle y \rangle x$. Using the same arguments as before, we construct other $3n_2 - 3$ terms in $H$, each one being a product of $n_2$ terms of $S$. It contradicts Eq. \eqref{sequencehi}, and implies that $S$ is not $2n$-product-one free.

\item Case $2\leq |S_{xy}|\leq 2n_2-1$. In this case, $|S_y| \geq 3n-2n_2+1>2n-1=\s(\langle y\rangle)$, therefore $S$ contains a $n$-product-one subsequence $R$ and $|S_y\bd R^{[-1]}| \geq 2n-2n_2+1$. If $S_y\bd R^{[-1]}$ is not $n$-product-one free, then we are done, thus we assume the contrary. Since $n_1>3$, we obtain that $k=2n_2-1\leq\lfloor n/2\rfloor+2$ and, as a consequence of Lemma \ref{lemmaegzinversecyclic}, $(y^a)^{[u]} \bd(y^b)^{[v]} \mid S_y$, where $\gcd(a-b,n)=1$ and $u,v\geq n-4n_2+5$. Since $n - 4n_2+5 \geq n_2$, we can assume, without lost of generality, that 
$$\begin{array}{ll}
h_i=y^{an_2}&\mbox{ if }i\in[1,2n_1-1];\\
h_i=y^{bn_2}&\mbox{ if }i\in[2n_1,3n_1-2];\\
h_i=y^{r_1n_2}x&\mbox{ if }i=3n_1-1.
\end{array}$$ 
The latter implies that $|T_{3n_1-1}\cap\langle y\rangle x| \ge 1$. Suppose that $e_i \in \langle y\rangle$ for every $i\in[1,n_2]$. If $|T_{3n_1-1}\cap\langle y\rangle x|=1$, then $S$ is not $2n$-product-one free, since the $2n_2-1$ terms in $\langle y \rangle$ generate a $n_2$-product-one subsequence over $\overline{G_{n,s}}$. Moreover, if $|T_{3n_1-1} \cap \langle y\rangle x| \geq 2$, say $g_{3n_1-1,1}, g_{3n_1-1,2} \in \langle y \rangle x$, then $\ldotprod_{1\leq i\leq n_2}e_i\bd\ldotprod_{2\leq j\leq n_2}g_{3n_1-1,j}$ has a subsequence $T_0$ of length $n_2$ whose product in some order belongs to $\langle y^{n_2} \rangle x$ and at least one of the remaining elements $e'_j \mid \left( \ldotprod_{1 \le i \le n_2} e_i \bd T_{3n_1-1} \bd T_0^{[-1]} \right)$ belongs to $\langle y\rangle x$. Therefore we can assume that $e_1 = y^rx$. We undo the subsequences $T_1 = (y^a)^{[n_2]}$ and $T_{2n_1} = (y^b)^{[n_2]}$ in order to create another subsequence of length $n_2$ from $(y^a)^{[n_2-1]} \bd (y^b)^{[n_2-1]} \bd (y^rx)$ whose product belongs to $D_{2n_1}$. Since $\gcd(a - b, n) = 1$, such subsequence attains a product in $\langle y^{n_2} \rangle x$. The other terms can be grouped in subsequences of length $n_2$ whose products belong to $D_{2n_1}$, and $n_2$ of them are left out. These $3n_1-1$ new terms over $D_{2n_1}$ (each of them being a product of $n_2$ terms from $S$) contain at least two in $\langle y^{n_2} \rangle x$, therefore it avoids the form given by Lemma \ref{lemmaegzinversedihedral}. We conclude that $S$ is not $2n$-product-one free.

\item Case $2n_2\leq |S_{xy}|\leq 3n_2-2$. Since $|S_{xy}| \geq 2n_2$, at least two of the $h_i$'s are products containing at least one term in $\langle y \rangle x$ each, therefore we may assume that $h_i=1$ for some $i \in [1,3n_1-2]$. It implies that either $h_i = 1$ for every $i \in [1,2n_1-1]$ or $h_i = 1$ for every $i \in [2n_1,3n_1-2]$. Consider the subcases:

\begin{enumerate}[{(3.}1)]
\item Subcase $|S_y\bd S_{y^{n_1}}^{[-1]}|\geq n_2-1$. Notice that the terms from $S_y\bd S_{y^{n_1}}^{[-1]}$ do not commute with the terms from $S_{xy}$. Since $|S_y \bd S_{y^{n_1}}^{[-1]} \bd S_{xy}| \geq 3n_2-1$, there exist two $n_2$-product-one subsequences over $\overline{G_{n,s}}$ formed by these terms. Therefore, some of the products $h_i$'s with $i \in [1,3n_1-2]$ can be reordered obtaining a distinct product, which contradicts Eq. \eqref{sequencehi} and then $S$ is not $2n$-product-one free.

\item Subcase $|S_y \bd S_{y^{n_1}}^{[-1]}| \leq n_2-2$. Suppose that $h_i=1$ for $i \in [1,2n_1-1]$ (resp. $i\in[2n_1,3n_1-2]$) and let $A = T_{2n_1} \bd \dots \bd T_{3n_1-1} \bd \ldotprod_{1 \leq i \leq n_2} e_i$ (resp. $A = T_{1} \bd \dots \bd T_{2n_1-1} \bd T_{3n_1-1} \bd \ldotprod_{1 \leq i \leq n_2} e_i$). Therefore, we have that $T_j \cap \langle y \rangle x = \varnothing$ for $j \in [2n_1,3n_1-2]$ (resp. $j\in[1,2n_1-1]$), otherwise we would have $h_i = 1$ for every $i \in [1,3n_1-2]$, a contradiction. 
Hence, $|A \cap \langle y\rangle x| \le 2n_2$. Then, since $|A|=n+n_2$ (resp. $|A|=2n+n_2$), it follows that $|A \cap \langle y \rangle| \geq n-n_2+1$ (resp. $|A\cap\langle y\rangle| \geq 2n-n_2+1$). Since $|S_y\bd S_{y^{n_1}}^{[-1]}|\leq n_2-2$, we obtain $|A \cap \langle y^{n_1} \rangle| \geq n-2n_2+2 = (n_1-2)n_2+2$ (resp. $|A\cap\langle y^{n_1}\rangle| \geq 2n-2n_2+2 = (2n_1-2)n_2+2$). Since $n_1 \geq 5$, it follows that $|A \cap \langle y^{n_1}\rangle| \geq 3n_2+2$ (resp. $|A \cap \langle y^{n_1}\rangle|\geq 8n_2+2$), therefore $A$ contains a subsequence of length $n_2$ whose product belongs to $H \cap \langle y^{n_1} \rangle = \{1\}$. Finally, we conclude that $S$ is not $2n$-product-one free. 
\end{enumerate}
\end{enumerate}
\qed

\subsection{The inverse problem: Proof of Theorem \ref{thmesinverse}(c)}\label{subsectioninversenodd}

Let $S \in \F(G_{n,s})$ be a $2n$-product-one free sequence of length $|S| = 3n-1$, and let $H = \langle x,y^{n_2} \rangle$. Similar to the Subsection \ref{subsectiondirectnodd}, we obtain that  $S = T_1 \bd \cdots \bd T_{3n_1-1} \bd E$, where 
$${T}_i={g}_{i,1}{\bd}\cdots{\bd}{g}_{i,n_2}, \quad {h}_i={g}_{i,1}\cdots{g}_{i,n_2}\in \pi(T_i) \subset H, \quad E = e_1 \bd \cdots\bd e_{n_2-1}$$ 
and 
\begin{equation}\label{eq:hs}
{h}_1\bd\cdots\bd{h}_{3n_1-1}=(y^{an_2})^{[2n_1-1]}\bd ({y}^{bn_2})^{[n_1-1]}\bd{y}^{r_1n_2}{x} \in \F(\overline{G_{n,s}}),
\end{equation} 
with $\gcd (a-b,n_1)=1$. First, we are going to deal with the extremal $2n$-product-one free sequences, and then with the other cases.

Suppose that $|S_{xy}|=1$, that is, $|S_y|=3n-2$. Since $\s(\langle y\rangle)=2n-1$, $S_y$ contains a $n$-product-one subsequence $R$. Since $S$ is $2n$-product-one free, $S_y\bd R^{[-1]}$ must be $n$-product-one free. By Lemma \ref{lemmaegzinversedihedral}, $S_y \bd R^{[-1]} = (y^a)^{[n-1]} \bd (y^b)^{[n-1]}$, where $\gcd(a-b,n)=1$. It implies that $(y^a)^{[n-1]} \bd (y^b)^{[n-1]} \mid S_y$. If there exists $g \mid S_y$ such that $g \not\in \{y^a, y^b\}$, then the sequence $S_y \bd (g \bd y^a \bd y^b)^{[-1]}$ contains a $n$-product-one subsequence $R'$. Notice that $|S_y \bd R'^{[-1]}| = 2n-2$, therefore we must have $S_y \bd R'^{[-1]} = (y^a)^{[n-1]} \bd (y^a)^{[n-1]}$ in order to $S_y \bd R'^{[-1]}$ be $n$-product-one free, but this is a contradiction since $g \mid S_y \bd R'^{[-1]}$. Thus, $S_y = (y^a)^{[u]} \bd (y^b)^{[v]}$, where $\gcd(a-b,n)=1$, and $u,v \geq n-1$. It implies that 
$$S=(y^a)^{[2n-1]}\bd(y^b)^{[n-1]}\bd y^rx,$$ 
and it is easy to verify that the latter is $2n$-product-one free.

If $|S_{xy}| = 0$, then $S \in \F(\langle y \rangle)$. Since $\s(\langle y \rangle) = 2n-1$, it is possible to extract two disjoint $n$-product-one subsequences, $S$ is not $2n$-product-one free. If $|S_{xy}| \ge 3n_2-1$, then we use the same argument than the case (1) of previous subsection in order to prove that $S$ is not $2n$-product-one free. If $|S_{xy}| = 2$, then necessarily $e_i \in \langle y \rangle x$ for some $i \in [1,n_2-1]$. Therefore the proof follows the same steps than case (2) of previous subsection and $S$ is not $2n$-product-one free.


From now on, we assume that $3 \leq |S_{xy}| \leq 3n_2-2$. Since $|S_{xy}| \geq 3$, it is possible to reorganize the subsequences $T_i$ in such way that $| (T_{3n_1-1} \bd E) \cap \langle y\rangle x | \geq 3$. Without loss of generality, we may assume that $g_{3n_1-1,1} \in \langle y \rangle x$. If $T_{3n_1-1} \bd g_{3n_1-1,1}^{[-1]} \bd E$ has a $n_2$-product-one subsequence over $\overline{G_{n,s}}$, then we reindex if needed and assume that $e_1\in\langle y\rangle x$. In this way, the proof that $S$ is not $2n$-product-one free follows the same steps than cases (2) and (3) of Subsection \ref{subsectiondirectnodd}. Therefore, suppose that $T_{3n_1-1} \bd g_{3n_1-1,1}^{[-1]} \bd E$ is $n_2$-product-one free over $\overline{G_{n,s}}$, and furthermore, $e_i \in \langle y\rangle$ for every $i\in[1,n_2-1]$. The same argument applies if $g_{3n_1-1,1}$ is replaced by any $g_{3n_1-1,i} \in \langle y \rangle x$. As a consequence of Lemma \ref{lemmaegzinversecyclic}, it follows that
$$g_{3n_1-1, i} = x^{\varepsilon_i} y^{a+n_2 \alpha_i} \;\;\;\;\; \mbox{ and } \;\;\;\;\; e_j = y^{b+n_2\beta_j},$$ 
where $i \in [1,n_2]$, $j \in [1,n_2-1]$, $\varepsilon_i \in \{0,1\}$, $\alpha_i, \beta_j \in [0, n_1-1]$, and $a,b \in [0, n_2-1]$ with $\gcd(a-b, n_2) = 1$. 

We claim that for every $k\in[1,3n_1-2]$, 
$$g_{k,i}=x^{\varepsilon_{k,i}}y^{a_k+n_2\alpha_{k,i}}$$ 
for every $i \in [1,n_2]$. In fact, if $T_k \bd g_{k,1}^{[-1]} \bd E$ is $n_2$-product-one free over $\overline{G_{n,s}}$ for some $k \in [1,3n_1-2]$, then claim holds by Lemma \ref{lemmaegzinversecyclic}. On the other hand, if $T_k \bd g_{k,1}^{[-1]} \bd E$ is not $n_2$-product-one free; say $T' \mid T_k \bd g_{k,1}^{[-1]} \bd E$ satisfies $\pi(T') \subset H$. Let $E' = T_k \bd E \bd T'^{[-1]}$ and $e_i' \mid E'$. The argument from the previous paragraph applies to show that the terms from $E'$ belong to the same class (that is, the exponents of $y$ belong to the same class modulo $n_2$), therefore either there exists $i_0 \in [2,n_2]$ such that $g_{k,i_0} = x^{\varepsilon_{k,i_0}} y^{b+n_2\alpha_{k,i_0}}$ and $g_{k,j} = x^{\varepsilon_{k,j}} y^{a_k+n_2\alpha_{k,j}}$ for $j \in [1,n_2] \backslash \{i_0\}$, or $g_{k,i} = x^{\varepsilon_{k,i}} y^{b+n_2\alpha_{k,i}}$ for every $i \in [1,n_2]$.

Since $h_1\neq h_{2n_1}$, we obtain that $a_1\not\equiv a_{2n_1} \pmod {n_2}$. Since $S$ is $2n$-product-one free, it follows that either $a_1\equiv a \pmod {n_2}$, or $a_{2n_1}\equiv a \pmod {n_2}$. Consider $R_0 = S \bd g_{1,1}^{[-1]} \bd g_{2n_1,1}^{[-1]} \bd g_{3n_1-1,1}^{[-1]}$. By the argument of the beginning of this subsection, we can write 
$$R_0 = R_1 \bd \dots \bd R_{3n_1-2} \bd \ldotprod_{1\le i \le 2n_2-4}\ell_i,$$ 
where $\pi(R_j) \subset H$ and $|R_j| = n_2$ for every $j \in [1, 3n_1-2]$, and the $\ell_i$'s are terms of $S$. 
Consider the subsequence $g_{1,1} \bd g_{2n_1,1} \bd g_{3n_1-1,1} \bd \ldotprod_{1\le i \le 2n_2-4} \ell_i$. The same algorithm can be applied for the sequences $S \bd g_{i,1}^{[-1]} \bd g_{j,1}^{[-1]} \bd g_{3n_1-1,1}^{[-1]}$ for every $i\in[1,2n_1-1]$, $j\in[2n_1,3n_1-2]$, obtaining either 
\begin{align*}
a_1\equiv\cdots\equiv a_{2n_1-1} &\equiv a \pmod {n_2}\;\;\;\;\;\mbox{ and }\;\;\;\;\; a_{2n_1}\equiv\cdots\equiv a_{3n_1-2}\equiv b \pmod {n_2}, \quad \text{ or } \\
a_1\equiv\cdots\equiv a_{2n_1-1} &\equiv b \pmod {n_2}\;\;\;\;\;\mbox{ and }\;\;\;\;\; a_{2n_1}\equiv\cdots\equiv a_{3n_1-2}\equiv a \pmod {n_2}.
\end{align*}

Assume that $a_1\equiv\cdots\equiv a_{2n_1-1}\equiv a \pmod {n_2}$; the other case is similar. Without loss of generality, suppose that $g_{1,1}, g_{2,1} \in \langle y\rangle$ (it is possible since $g_{i,1} \dots g_{i,n_2} \in \langle y^{n_2} \rangle$ for $i \in \{1,2\}$, that is, the terms $g_{i,j} \in \langle y \rangle x$ cancel the $x$ in pairs, and $n_2$ is odd). Since $E \in \F(\langle y \rangle)$ and $|(T_{3n_1-1} \bd E) \cap \langle y \rangle x| \geq 3$, it is possible to select (say) $g_{3n_1-1,1} \bd g_{3n_1-1,2} \mid (T_{3n_1-1} \cap \langle y \rangle x)$. Since $a_1\equiv a_2\equiv a \pmod {n_2}$, we change some terms, and finally obtain the following three products in $H \cap \langle y \rangle x$ of length $n_2$ each:
$$g_{1,1} \cdot g_{2,1} \cdot \prod_{3 \leq i \leq n_2} g_{3n_1-1,i} , \quad g_{3n_1-1,1} \cdot \prod_{2 \leq i \leq n_2} g_{1,i} , \quad g_{3n_1-1,2} \cdot \prod_{2 \leq i \leq n_2} g_{2,i}.$$ 
It contradicts the form given by Eq. \eqref{eq:hs}, hence $S$ is not $2n$-product-one.
\qed

\end{document}